\documentclass[
  pre,
  aps,
  a4paper,
  english,
  showpacs,
  reprint,
  twocolumn,
  superscriptaddress
]{revtex4-1}

\usepackage[T1]{fontenc}
\usepackage[utf8]{inputenc}
\usepackage{amsmath}
\usepackage{amssymb}
\usepackage{graphicx}
\usepackage{refstyle}
\usepackage[
  citecolor=blue,
  colorlinks,
  linkcolor=blue,
  urlcolor=blue,
]{hyperref}
\usepackage[dvipsnames]{xcolor}

\begin{document}

\title{Conservative Perturbation Theory for Non-Conservative Systems}
\author{Tirth Shah}
\email{ph12b007@smail.iitm.ac.in}
\affiliation{
  Department of Physics,
  Indian Institute of Technology Madras,
  Chennai, Tamilnadu 600036, India
}
\author{Rohitashwa Chattopadhyay}
\email{crohit@iitk.ac.in}
\affiliation{
  Department of Physics,
  Indian Institute of Technology Kanpur,
  Uttar Pradesh 208016, India
}
\author{Kedar Vaidya}
\email{kedarv93@vt.edu}
\affiliation{
Department of Biomedical Engineering and Mechanics,
Virginia Tech,
Blacksburg, VA 24061, USA
}
\author{Sagar Chakraborty}
\email{sagarc@iitk.ac.in}
\affiliation{
  Department of Physics,
  Indian Institute of Technology Kanpur,
  Uttar Pradesh 208016, India
}
\affiliation{
  Mechanics and Applied Mathematics Group,
  Indian Institute of Technology Kanpur,
  Uttar Pradesh 208016, India
}
%

%%% Abstract -----------------------------------------------------------

\begin{abstract}
In this paper, we show how to use canonical perturbation theory for dissipative dynamical systems capable of showing limit cycle oscillations.
Thus, our work surmounts the hitherto perceived barrier for canonical perturbation theory that it can be applied only to a class of conservative systems, viz.,~Hamiltonian systems.
In the process, we also find Hamiltonian structure for an important subset of Li\'enard system--- a paradigmatic system for modeling isolated and asymptotic oscillatory state.
We discuss the possibility of extending our method to encompass even wider range of non-conservative systems.
\end{abstract}
%%%%%%%%%%%%%%
%%%%%%%%%%%%%%%%%
\maketitle
%%%%%%%%%%%%%%%%%%%%%

\section{Introduction}
Limit cycle~\cite{JEN2013}, an intriguing mathematical entity, is the subject of the still-unsolved second part of celebrated Hilbert's sixteenth problem~\cite{ILY02}.
The initial driving question behind this paper was: can one quantize a dynamical system that possesses a stable (isolated) limit cycle in its phase space.
In other words, how would the system which for a set of initial conditions always asymptotically reaches a fixed oscillatory state, quantum mechanically behave?
Many systems, e.g.,~simple harmonic oscillator, have non-isolated family of uncountably many closed phase orbits around center-type fixed points.
This in simpler language means that these systems can oscillate with amplitudes that can be continuously varied by continuously changing the initial conditions.
It is now a common knowledge that such systems can be quantized leading to quantized oscillatory states labeled by discrete quantum numbers. 
However, it is not at all clear what an attempt to quantize a system with countable number of isolated closed phase trajectories (limit cycles) will get us.
A bit of thought straightaway suggests that any such attempt will be infested with many fundamental obstacles: (i) the systems having stable limit cycles are \emph{nonlinear}, (ii) such systems are \emph{dissipative} (phase-volume-contracting), (iii) these systems are \emph{non-Hamiltonian} in general, etc.

In view of the above, we set a comparatively modest but non-trivial goal in this paper which potentially is the first stepping stone towards our aforementioned greater aim.
\emph{Our goal is to formulate a Hamiltonian formalism for systems with limit cyclic behaviour.} 
One of the immediate probable benefit may be that one can tap into the machinery of canonical perturbation theory (used exclusively for conservative systems having a Hamiltonian description) and explore dissipative systems with its help.
In this manner, one can contribute to the plethora of useful basic perturbative techniques available for such systems~\cite{AHN2000}.
It is worthwhile to remark that the invaluable canonical perturbation theory (cf.~\cite{JE15} for a very recent application) is a cornerstone in the centuries-old subject of classical mechanics. 
The perturbative technique is an indispensable analytical tool since most of the realistic problems--- classical three-body problem~\cite{CEN07} being one of the most celebrated one--- are not exactly solvable.
We showcase in this paper the fact that we have actually obtained a Hamiltonian formalism for an important class of nonlinear dissipative systems.
Thence, \emph{we have invented an apparently oxymoronic perturbative technique: canonical perturbation method for dissipative systems.}
This is the primary result of this paper.

In order to pursue our redefined goal systematically and for the sake of concreteness, we have chosen to illustrate our main idea and results using the widely studied system possessing a limit cycle, viz.,~Van der Pol oscillator (VdPO).
VdPO and its modified forms are of practical utility in physical~\cite{PO93,HAY98,BEN13}, chemical~\cite{Lee13}, and biological sciences~\cite{EU06,GLA95}; in mathematics~\cite{Wec07}, and economics \cite{AK86}; and of late, even in seismology \cite{CEGP99}, and solar cycle~\cite{MIN2000}; not to mention the variety of engineering disciplines.
As shown later, the results obtained in this paper using VdPO can be straightforwardly extended to encompass a wide class of Li\'enard systems.
Our idea is motivated by the Hamiltonian description available for a damped simple harmonic oscillator (DSHO) which, however, is a linear system incapable of showing any limit cycle behaviour.

The fundamental concept of Hamiltonians for conservative systems with even dimensional phase space has voluminous literature at every possible level. 
Hamiltonian description for dissipative systems is rather rare for obvious reasons. 
Bateman's dual Hamiltonian~\cite{BAT1931,DEK1981,UM2002} is one such Hamiltonian tailor-made for DSHO. 
Bateman's dual Hamiltonian yields two second order differential equations of motion: one of the equations is the usual equation of motion for DSHO while the other one is an \emph{auxiliary} equation. 
Together they conserve a global time-independent first integral of motion: the Bateman's dual Hamiltonian.
Also, an equivalent Hamiltonian--- Caldirola--Kanai Hamiltonian--- for DSHO is obtained by replacing the constant mass in the Hamiltonian for simple harmonic oscillator with an exponentially varying mass.
Alternatively, among some other possibilities~\cite{HAV1957}, one can also use fractional derivatives to write a Hamiltonian for DSHO~\cite{RIE96}.

A somewhat different kind of Hamiltonian formalism to deal with non-conservative systems described by a set of ordinary differential equations has been developed~\cite{GAL2013}. 
The motivation for such a work was the apparent failure of the classical Rayleigh--Lagrange equations to account for those dissipative forces which are nonlinear in velocities. 
It has however been established that the nonlinear dissipative forces can be treated in the variational formalism by extending the Rayleigh--Lagrange framework~\cite{VIR2015}. 
Once the Hamiltonian of a system is obtained, the system can, in principle, be quantized by treating the observables as operators in an appropriate Hilbert space and invoking Schr\"odinger's equation.
In this spirit, the quantization of the damped harmonic oscillator is achieved by quantising the Bateman's dual Hamiltonian~\cite{FES1977}. 
This leads to justifiable~\cite{GUE2012,Chru2006} complex eigenvalues for the Hamiltonian. 
It has also been shown that the quantum DSHO is equivalent to the quantum problem with a two-dimensional parabolic potential barrier~\cite{CHR2006}.
Attempts have also been made to quantize DSHO via a Hamiltonian formalism different from the Bateman's dual Hamiltonian~\cite{UM2002,CHA2007}, e.g.,~the path integral approach, quantising the Hamiltonian obtained from the integrals of motion via the modified Prille-Singer approach etc.

In this paper, we tap into the non-quantum mechanical aspects of the Hamiltonian of DSHO to  develop the Hamiltonian description of VdPO and related non-conservative systems. 
While we plan to pursue the quantum-mechanical aspect of dissipative systems in detail in future, little thought reveals --- as we also show explicitly later in the paper--- that in the neighbourhood of the limit cycle, an effective Hamiltonian of VdPO can be obtained by taking the average of the Hamiltonian of VdPO over the limit cycle. 
The effective Hamiltonian will give a linearized equation of motion mimicking a DSHO.
\emph{Thus, we believe that the aforementioned quantum mechanics of DSHO naturally, although approximately, carries over to VdPO in the neighbourhood of the limit cycle.} 

We now proceed to present the Hamiltonian formalism for VdPO and lead the reader through its canonical perturbation theory to find the correction to the frequency upto second order in the strength of nonlinearity. 
But before that, for the sake of clarity and completeness, we shall go through a brief discussion of the Hamiltonian formalism of DSHO tuned to our purpose.
%
%%%%%%%%%%%%%%%%%%%%%%%%%%%%%%%%%%%%%%%%%%%%%%%%%%%%%%%%%%%%%%%%%%%%
%
\section{Hamiltonian formalism for DSHO}
 The equation of motion for DSHO of unit mass is given by,
\begin{equation}
\label{eqdsho}
\ddot{x}+2\lambda\dot{x}+\Omega^2{x}=0\,. 
\end{equation}
The strategy employed by Bateman was to add a negatively damped linear oscillator,
\begin{equation}
\quad \ddot{y}-2\lambda\dot{y}+\Omega^2{y}=0\,,\label{eqadsho}
\end{equation}
to it so that the total system becomes conservative and the conserved quantity, interestingly enough, acts as a Hamiltonian ($H_B$) of the system.
The Lagrangian for DSHO is given by,
\begin{equation}
L(x,y,\dot{x},\dot{y})= \dot{x}\dot{y}-\Omega^2 {x}{y}+\lambda(x\dot{y}-\dot{x}y).
\end{equation}
Using Legendre transformation, the Bateman's Dual Hamiltonian, is found out to be:
\begin{equation}
\label{BDH}
H_B(x,y,p_x,p_y)= p_{x}p_{y}+\omega^2 xy-\lambda(p_{x}x-p_{y}y)\,,
\end{equation}
where $p_x=\partial L/\partial \dot{x}$, $p_y=\partial L/\partial\dot{y}$, and $\omega^2\equiv(\Omega ^2-\lambda^2)$.
This Hamiltonian gives back both Eq.~(\ref{eqdsho}) and Eq.~(\ref{eqadsho}).

Following (time-dependent) Caldirola-Kanai Hamiltonian also gives equation of motion for DSHO:
\begin{equation}
\label{CK}
H_{CK}(x,p_x)=\frac{\exp(-2 \lambda t)}{2}p_x^2+\frac{\exp(2 \lambda t)}{2}\Omega^2 x^{2}.
\end{equation}
$H_{CK}$ can be interpreted as the Hamiltonian of a simple harmonic oscillator with mass changing exponentially with time.
It is equivalent to $H_B$ as it can be obtained from $H_B$ by first performing a series of canonical transformations: $\left(x,y,p_{x},p_{y}\right)\rightarrow\left(Q_1,Q_2,P_1,P_2\right)\rightarrow\left(Q_{11},Q_{22},
P_{11},P_{22}\right)$ effected respectively by the generating functions:
$F_{2a}(x,y,P_1,P_2)=xP_1 \exp(\lambda t)+yP_2\exp(-\lambda t)$ and $F_{2b}\left(Q_1,Q_2,P_{11},P_{22}\right)=P_{11}\left({Q_1+Q_2}\right)/\sqrt{2}+P_{22}\left({Q_1-Q_2}\right)/\sqrt{2}$; and then by using the generating functions 
$F_{2c}(Q_{11},\pi)=\pi Q_{11}+\frac{\lambda Q^2_{11}}{2}$ and
$F_{2d}\left(\xi,P\right)=P\xi \exp(-\lambda t)$ to respectively generate canonical transformations:
$(Q_{11},P_{11})\rightarrow (\xi,\Pi)\rightarrow(Q,P)$. 
This results in $H_{CK}(Q,P)$.

An immediate application of $H_B$ would be to employ Hamilton-Jacobi theory to solve DSHO completely for arbitrary initial conditions. It is quite straightforward as the Hamilton-Jacobi equation turns out to be separable on performing the sequence of canonical transformations using $F_{2a}$ and $F_{2b}$. 
\section{Hamiltonian formalism for V\lowercase{d}PO}
The well-known Van der Pol oscillator is governed by the following equation:
\begin{equation}
\ddot{x}+\varepsilon(x^2-1)\dot{x}+\omega^2x=0\,.\label{vdpo}
\end{equation}
This two dimensional autonomous dynamical system undergoes non-generic Hopf bifurcation at $(x,\dot{x},\varepsilon)=(0,0,0)$: stable limit cycle appears for positive $\varepsilon$, whereas an unstable limit cycle appears for negative $\varepsilon$.
Existence of attractors--- either stable focus or stable limit cycle--- in the phase space rules out any Hamiltonian for VdPO.
Our strategy is, thus, to introduce a two dimensional autonomous auxiliary dynamical system complementing VdPO in such a fashion that a Hamiltonian, \`a la Bateman's dual Hamiltonian for DSHO, may be written down for the resulting complete four-dimensional system.
\subsection {The Hamiltonian}
Since searching for the Hamiltonian is effectively an exercise in trial and error at this juncture, it helps to start by considering the following set of symmetrical (except for the sign in front of the second terms) equations:
\begin{subequations}
\begin{eqnarray}
&&\ddot{x}+\varepsilon(\alpha x^2+\beta y^2-1)\dot{x}+\omega^2x=0\,,\\
&&\ddot{y}-\varepsilon(\alpha x^2+\beta y^2-1)\dot{y}+\omega^2y=0\,.
\end{eqnarray}
\label{ab}
\end{subequations}
Here $\alpha$ and $\beta$ are real parameters.
One can see that this is a conservative system by noting that the divergence of the vector field $(\dot{x},\ddot{x},\dot{y},\ddot{y})$ is zero. 
Therefore, any four dimensional hyper-volume of initial conditions keeps its measure intact as it evolves in the corresponding phase space.
This raises hope of finding a Hamiltonian formulation for it.
In fact, the function,
\begin{eqnarray}
&&L(x,y,\dot{x},\dot{y})= \dot{x}\dot{y}-\omega^2xy+\nonumber\\
&&\phantom{L}\frac{\varepsilon}{2}\left[\left(y-\alpha yx^2-\beta \frac{y^3}{3}\right)\dot{x}\right.- \left.\left(x-\beta xy^2-\alpha\frac{x^3}{3}\right)\dot{y}\right]\quad\,
\end{eqnarray}
can be seen to yield Eqs. (\ref{ab}) via Euler-Lagrange equations, and hence can be consider as a valid Lagrangian for the system under consideration.

Putting $\beta=0$ in Eqs. (\ref{ab}), we get --- along with an auxiliary equation --- Van-der-Pol oscillator (VdPO) in $x$:
\begin{subequations}
\begin{eqnarray}
&&\ddot{x}+\varepsilon(\alpha x^2-1)\dot{x}+\omega^2x=0\,,\label{VdPO}\\
&&\ddot{y}-\varepsilon(\alpha x^2-1)\dot{y}+\omega^2y=0\,.\label{auxiliary VdPO}
\end{eqnarray}
\end{subequations}
We call the immediately preceding driven equation in $y$ as auxiliary VdPO. 
It is important and insightful to realize that these two equations are coupled \emph{unidirectionally}: VdPO, which is of our prime interest, is not affected by the dynamics of auxiliary VdPO.
We haven't put the bookkeeping parameter $\alpha=1$ just to keep track of the nonlinear terms. 
The appropriate Lagrangian ($L_v$, say) is then
\begin{eqnarray}
L_v=\dot{x}\dot{y}-\omega^2xy+
\frac{\varepsilon}{2}\left[\left(y-\alpha yx^2\right)\dot{x}- \left(x-\alpha\frac{x^3}{3}\right)\dot{y}\right].\qquad
\end{eqnarray}
The conjugate momenta are easily obtained as follows:
\begin{subequations}
\begin{eqnarray}
&&p_x=\frac{\partial L}{\partial \dot{x}}=\dot{y}+\frac{\varepsilon}{2}\left(y-\alpha yx^2\right)\,,\\
&&p_y=\frac{\partial L}{\partial \dot{y}}=\dot{x}-\frac{\varepsilon}{2}\left(x-\alpha\frac{x^3}{3}\right)\,.
\end{eqnarray}
\end{subequations}
Hence, the corresponding Hamiltonian, $H_v\,(\textrm{say})=\dot{x}p_x+\dot{y}p_y-L_v$ (standard Legendre transformation), is
\begin{eqnarray}
H_v&=&\left[p_xp_y+\frac{\varepsilon}{2}\left(xp_x-yp_y\right)+\left(\omega^2-\frac{\varepsilon^2}{4}\right)xy\right]
\nonumber \\ 
&&+\alpha\left[\frac{\varepsilon}{2}\left(x^2yp_y-\frac{x^3p_x}{3}\right)+\frac{\varepsilon^2}{12}\left(4x^3y-\alpha x^5y\right)\right]\,.\qquad\quad\label{hvdpo}
\end{eqnarray}
Using $H_v$, one can directly verify that the standard canonical equations of motion yield VdPO and auxiliary VdPO.
\subsection{Dynamics of 4D VdPO--auxiliary VdPO system}
In principle, one can use all the conventional machinery of analysing nonlinear dynamical systems to study the dynamics of four dimensional autonomous VdPO--auxiliary VdPO system.
Even though that is not our main interest, it is quite interesting to investigate this aspect as well.

For the sake of convenience if we choose $\alpha=1$, Eq.~(\ref{VdPO}) allows for an asymptotic solution (which is stable if $\varepsilon>0$): $x\approx2\cos\omega t$ as $|\varepsilon|\rightarrow0$. 
Since the two equations are unidirectionally coupled, the large-time solution of Eq.~(\ref{auxiliary VdPO}) should be approximately governed by 
\begin{eqnarray}
&&\ddot{y}+\varepsilon\left(1-4\cos^2\omega t\right)\dot{y}+\omega^2y=0\,,
\end{eqnarray}
wherein substituting 
\begin{eqnarray}
y&&=u\exp\left[-\varepsilon\int \left(\frac{1}{2}-2\cos^2\omega t\right) dt\right]\nonumber \\&&=u\exp\left[\frac{\varepsilon}{2} \left(t+\frac{\sin2\omega t}{\omega}\right)\right]\,,\label{subs}
\end{eqnarray}
one obtains following Hill's equation
\begin{eqnarray}
\ddot{u}+p(t)u=0\,,
\end{eqnarray}
where 
\begin{eqnarray}
&&p(t)=\omega^2-\frac{\varepsilon^2\left(1-4\cos^2\omega t\right)^2}{4}-{2\varepsilon\omega\sin2\omega t}\,,\quad\\
\Rightarrow\quad&&p(t)\approx\omega^2-{2\varepsilon\omega\sin2\omega t}+\mathcal{O}(\varepsilon^2)\,,
\end{eqnarray}
is a periodic function: $p(t+2\pi/\omega)=p(t)$.
The preceding approximated form for $p(t)$ suggests strong parametric resonance in the corresponding Mathieu equation in $u$ which through Eq.~(\ref{subs}) asserts that the solution for $y$ grows rapidly for even small enough positive $\varepsilon$.

However, the unidirectional nature of the coupling between VdPO and auxiliary VdPO allows us to ignore this resonance since VdPO remains unaffected by it.
\subsection{Equivalent linearization}
We know that there exists a limit cycle solution for Eq.~(\ref{VdPO}). 
Now, let us seek the equivalently linearized solution in the neighbourhood of the limit cycle.
To this end, we write an equivalent quadratic Hamiltonian ($H_{2v}$, say) starting from $H_v$:
\begin{eqnarray}
H_{2v}=\left[p_xp_y+\frac{\varepsilon}{2}\left(xp_x-yp_y\right)+\left(\omega^2-\frac{\varepsilon^2}{4}\right)xy\right]\nonumber \qquad\qquad\quad\\
+\alpha\left[\frac{\varepsilon}{2}\left(\langle x^2\rangle yp_y-\frac{\langle x^3\rangle p_x}{3}\right)+\frac{\varepsilon^2}{12}\left(4\langle x^3\rangle y-\alpha\langle  x^5\rangle y\right)\right],\qquad
\end{eqnarray}
where $\langle\cdots\rangle$ is the average over one period of the limit cycle $x\sim A\cos\omega t$. 
Hence,
\begin{eqnarray}
H_{2v}=&&p_xp_y+\frac{\varepsilon}{2}\left\{xp_x-\left(1-\frac{\alpha A^2}{2}\right)yp_y\right\}\nonumber \\ 
&&+\left(\omega^2-\frac{\varepsilon^2}{4}\right)xy .
\end{eqnarray}
This effective Hamiltonian gives the following linearized equation in $x$
\begin{eqnarray}
\ddot{x}+\varepsilon\left(\frac{\alpha A^2}{4}-1\right)\dot{x}+\omega^2x=0+\mathcal{O}(\varepsilon^2)\,,
\end{eqnarray}
which is the correct expression~\cite{JS07}. 
The accompanying equation for $y$ is
\begin{eqnarray}
\ddot{y}-\varepsilon\left(\frac{\alpha A^2}{4}-1\right)\dot{y}+\omega^2y=0+\mathcal{O}(\varepsilon^2)\,.
\end{eqnarray}
Thus, as expected in the neighbourhood of the limit cycle VdPO-auxiliary VdPO system becomes the standard system addressed by Bateman's dual Hamiltonian.
Hence, all the studies, whether classical or quantum, present in the literature on DSHO via Hamiltonian formalism naturally carries over to VdPO in the neighbourhood of the limit cycle.
\subsection{Generalizations}
In this subsection, we intend to discuss some straightforward but useful generalizations of the aforementioned results.
\subsubsection{Forced Van der Pol Oscillator}
Consider VdPO which is both externally and parametrically forced.
Without any loss of generality, we take the forces to be sinusoidal.
Therefore, we have
\begin{eqnarray}
\ddot{x}+\varepsilon(\alpha x^2-1)\dot{x}+\left(\omega^2+F_1\cos\gamma t\right)x=F_2\cos\Omega t\,,\quad\label{fvdpo}
\end{eqnarray}
In order to propose a Hamiltonian for it, we first guess following auxiliary equation
\begin{eqnarray}
\ddot{y}-\varepsilon(\alpha x^2-1)\dot{y}+\left(\omega^2+F_1\cos\gamma t\right)y=F_2\cos\Omega t\,.\quad\label{fauxiliary VdPO}
\end{eqnarray}
Now, by trivial inspection we note that a time-dependent function,   
\begin{eqnarray}
H_f(x,y,\dot{x},\dot{y},t)\equiv H_v-(x+y)F_2\cos\Omega t+xyF_1\cos\gamma t,\,\qquad
\end{eqnarray}
can be written down for the unidirectionally coupled 4D non-autonomous system under consideration.
\subsubsection{Li\'enard Systems}
The entire preceding exercise can be extended for a more general type of equations called Li\'enard equation~\cite{JS07}: 
\begin{eqnarray}
&&\ddot{x}+\varepsilon f(x)\dot{x}+\omega^2x=0\,,\label{le}
\end{eqnarray}
where $f(x)$ is a continuously differentiable function.
By virtue of Li\'enard's theorem, such systems may possess at least one periodic solution when conditions of the theorem are met.
VdPO is just a special case of  such a system.
One can verify that the Lagrangian
\begin{eqnarray}
\label{lag}
L=\left[\dot{x}\dot{y}+\frac{\varepsilon}{2}\left(\dot{x}y-x\dot{y}\right)-\omega^2xy\right]+{\varepsilon}\left[f_1(x)\dot{y}-f_2(x)\dot{x}y\right],\nonumber\\
\end{eqnarray}
corresponds to Eq.~(\ref{le})
along with a corresponding auxiliary differential equation in $y$. 
Here, one must recognize $f(x)$ as $f'_1(x)+f_2(x)-1$. 
Therefore, we note that in order to get Eq.~(\ref{le}), one can in fact work with different Lagrangians by only demanding that the combination of $f_1(x)$ and $f_2(x)$ gives back the desired $f(x)$.
This is merely the consequence of gauge invariance of Lagrangian: a Lagrangian is non-unique up to addition of a total time derivate of a function of $x$ and $y$.
The corresponding Hamiltonian is given by:
\begin{eqnarray}
H=p_{x}p_{y}+\frac{\varepsilon }{2}\left(p_{x}x-p_{y}y\right)-\varepsilon f_{1}\left(x\right)p_{x}+\varepsilon p_{y}yf_{2}\left(x\right) \nonumber\\-\frac{\varepsilon^2 xy}{4}+\frac{\varepsilon^2 f_1(x)y}{2}-\varepsilon^2 f_2(x)y+\frac{\varepsilon^2 f_2(x)f_1(x)xy}{2}+\omega^2xy.\nonumber\\
\end{eqnarray}
It goes without saying that forced Li\'enard equation naturally comes under the radar of Hamiltonian formalism along similar lines. 
Moreover, the more general Li\'enard equation that has $\omega^2x$ replaced by a continuous function of $x$, say $g(x)$, can also be analogously studied, specially when $g(x)$ can be written down as a derivative of some potential function.
%%%%%%%%%%%%%%%%%%%%%%%%%%%%%%%%%%%%%%%%%%%%%%%%%%%%%%%%%%%%%%%%%%%%
\section{Interpreting auxiliary equation}
Up to now we have made use of the auxiliary equation in a manner so as to be able to write down a Hamiltonian for the complete coupled system.
We have made it a point to keep the coupling between the system of interest and the auxiliary system unidirectional so that the dynamics of the former gets no feedback from that of the latter.
We have chosen to remain completely disinterested in the form or the properties of the auxiliary equation, so much so that \emph{in practice} it is completely unnecessary for us to even bother about the interpretation of the variables of the auxiliary equation.
Nonetheless, we can give a meaning to the auxiliary variables by connecting with Galley's modified formulation of Hamilton's principle: Hamilton's principle with initial data~\citep{GAL2013}.

To this end we define
\begin{eqnarray}
\label{q1q2}
q_{1}\equiv\frac{2x+y}{2},\, 
q_{2}\equiv\frac{2x-y}{2}\,.
\end{eqnarray}
We exploit the gauge invariance of Lagrangian in Eq.~(\ref{lag}) to redefine $L$ as
\begin{eqnarray}
&&L\rightarrow L+\frac{d}{dt}\left\{\frac{\varepsilon}{2}xy-\varepsilon f_1(x)y\right\}\,,\\
\Rightarrow\, &&{L}=\dot x\dot y-\omega^2 xy-\varepsilon f\left(x\right)\dot{x}y\,.
\end{eqnarray}
Corresponding Hamiltonian is
\begin{eqnarray}
H=p_{x}p_{y}+\omega^2 xy+\varepsilon f\left(x\right)yp_{y}\,,\label{cvdpo}
\end{eqnarray}
with conjugate momenta being: 
$p_{x}=\dot{y}-\varepsilon f(x)y$ and  $p_{y}=\dot{x}$. 
Subsequent use of definitions (\ref{q1q2}) in the expression for ${L}$ results in
\begin{eqnarray}
{L}=\left(\frac{ \dot{q}_{1}^2}{2}-\omega^2\frac{q_{1}^2}{2}\right)-\left(\frac{ \dot{q}_{2}^2}{2}-\omega^2\frac{q_{2}^2}{2}\right)+N\,.\label{lt}
\end{eqnarray}
Here, ${N}=N(q_1,q_2,\dot{q}_1,\dot{q}_2)$ is given by:
\begin{eqnarray}
N\equiv-\varepsilon \left[q_1-q_2\right]\left[\frac{\dot{q}_1+\dot{q}_2}{2}\right]f\left(\frac{q_1+q_2}{2}\right)\,.
\end{eqnarray}
We note that the two terms in the curly brackets in Eq.~(\ref{lt}) correspond to the conservative part of the Lagrangian.
Thus, in accordance with Hamilton's principle with initial data, the auxiliary variable (e.g.,~$y$ of auxiliary VdPO) combines with the variable of Li\'enard equation (e.g.,~$x$ of VdPO) to identify the function $N$.
$N$ describes the nonconservative forces and couples the two oppositely traversed stationary paths in the configuration space of the system of interest (e.g.,~VdPO) obtained in the \emph{physical limit}: $q_1=q_2$.
%%%%%%%%%%%%%%%%%%%%%%%%%%%%%%%%%%%%%%%%%%%%%%%%%%%
\section{Canonical Perturbation Theory for V\lowercase{d}PO}
We now come to the heart of the paper.
Having found a Hamiltonian description for VdPO, we are naturally very optimistic about employing canonical perturbation theory to perturbatively find the frequency of the limit cycle oscillation of VdPO.
To this end, it is very convenient to work with the simpler Hamiltonian given in Eq.(\ref{cvdpo}).
In principle, the Hamiltonian of Eq.~(\ref{hvdpo}) should yield the same end result. 
Hence, our starting point is
\begin{equation}
\label{H}
H(x,y,p_x,p_y)=p_xp_y+\omega^2xy+\varepsilon(x^2-1)yp_y.
\end{equation}
We do a canonical transformation $\left(x,y,p_{x},p_{y}\right)\rightarrow\left(X,Y,P_X,P_Y\right)$ generated by $F_{2b}\left(x,y,P_{X},P_{Y}\right)=P_{X}\left({x+y}\right)/\sqrt{2}+P_{Y}\left({x-y}\right)/\sqrt{2}$.
The new Hamiltonian $H(X,Y,P_X,P_Y)$ is then
\begin{eqnarray}
\label{HQP}
H\left(X,Y,P_X,P_Y\right)=\frac{P_X^2}{2}+\omega^2\frac{X^2}{2}-\left(\frac{P_Y^2}{2}+\omega^2\frac{Y^2}{2}\right)\nonumber \\+\varepsilon\left(\frac{P_X-P_Y}{\sqrt{2}}\right)\left(\frac{X-Y}{\sqrt{2}}\right)\left[\frac{1}{2}\left(X+Y\right)^2-1\right]\,.\nonumber
\end{eqnarray}
The next step in canonical perturbation theory demands us to rewrite $H(X,Y,P_X,P_Y)$ in terms of the action-angle variables $(\phi^{(0)}_{1},\phi^{(0)}_{2}, I^{(0)}_{1},I^{(0)}_{2})$ of the unperturbed ($\varepsilon=0$) Hamiltonian. 
These variables are related to the old variables $(X,Y,P_X,P_Y)$ as follows 
\begin{eqnarray}
\label{AA1}
\begin{split}
X=\sqrt{\frac{2I^{(0)}_{1}}{\omega}}\sin\phi^{(0)}_{1}, 
\quad P_X=\sqrt{2\omega I^{(0)}_{1}}\cos\phi^{(0)}_{1},\quad 
\\Y=\sqrt{\frac{2I^{(0)}_{2}}{\omega}}\sin\phi^{(0)}_{2}, 
\quad P_Y=\sqrt{2\omega I^{(0)}_{2}}\cos\phi^{(0)}_{2}.\quad
\end{split}
\end{eqnarray}
The new Hamiltonian, $K(\phi^{(0)}_1,\phi^{(0)}_2,I^{(0)}_1,I^{(0)}_2)$ say, would be
\begin{eqnarray}
\label{K0}
&&K=K_0(I^{(0)}_{1},I^{(0)}_{2})+\varepsilon K_1(\phi^{(0)}_1,\phi^{(0)}_2,I^{(0)}_1,I^{(0)}_2)\,,\\
\Rightarrow\,&&K=I^{(0)}_{1}\omega-I^{(0)}_{2}\omega+\varepsilon K_1(\phi^{(0)}_1,\phi^{(0)}_2,I^{(0)}_1,I^{(0)}_2)\,.\quad
\end{eqnarray}
Let $(I_1,I_2,\phi_1,\phi_2)$ be the action-angle variables of $H$. 
The generating function for the transformation $(\phi^{(0)}_1,\phi^{(0)}_2,I^{(0)}_1,I^{(0)}_2)\rightarrow(\phi_1,\phi_2,I_1,I_2)$ is written as a near-identity transformation in the powers of $\varepsilon$:
\begin{eqnarray}
\label{S}
S(\phi^{(0)}_1,\phi^{(0)}_2,I_1,I_2)=\phi^{(0)}_1 I_1+\phi^{(0)}_{2}I_2
+\varepsilon S_1(\phi^{(0)}_1,\phi^{(0)}_2,I_1,I_2)\quad\nonumber\\+\varepsilon^2 S_2(\phi^{(0)}_1,\phi^{(0)}_2,I_1,I_2)+\mathcal{O}(\varepsilon^3).\qquad\qquad
\end{eqnarray}
Therefore, the Hamiltonian $E(I_1,I_2)$ ($H$ written in terms of its action-angle variable) can be obtained as
\begin{equation}
\label{totE}
E(I_1,I_2)=E_0(I_1,I_2)+\varepsilon E_1(I_1,I_2)+\varepsilon^2 E_2(I_1,I_2)+\mathcal{O}(\varepsilon^3)\,,
\end{equation}
where~\cite{JS98} 
\begin{subequations}
\begin{eqnarray}
\label{enecorr0}
&&E_0(I_1,I_2)=K_0(I_1,I_2)=I_1\omega-I_2\omega,\\
\label{enecorr1}
&&E_1(I_1,I_2)=K_1+\sum\limits_{i=1}^2\nu^{(0)}_{i}\frac{\partial S_1}{\partial\phi^{(0)}_{i}}=\langle K_1\rangle=0,\\
\label{enecorr2}
&&E_2(I)=\sum\limits_{i=1}^2\left(\frac{\partial K_1}{\partial I_{i}}\frac{\partial S_1}{\partial \phi^{(0)}_{i}}+\nu^{(0)}_{i}\frac{\partial S_2}{\partial \phi^{(0)}_{i}}\right)\nonumber \\
&\phantom{=}&+\frac{1}{2}\sum\limits_{i,j=1}^2\frac{\partial \nu^{(0)}_{i}}{\partial I_{j}}\frac{\partial S_1}{\partial \phi^{(0)}_{i}}\frac{\partial S_1}{\partial \phi^{(0)}_{j}}= \left\langle\sum\limits_{i=1}^2\frac{\partial K_1}{\partial I_{i}}\frac{\partial S_1}{\partial \phi^{(0)}_{i}}\right\rangle.\qquad\label{E_2}
\end{eqnarray}
\end{subequations}
Here, $(\nu^{(0)}_1,\nu^{(0)}_2)\equiv({\partial E_0}/{\partial I_{1}},{\partial E_0}/{\partial I_{2}})=(\omega,-\omega)$, $\left\langle \cdots \right\rangle$ is defined as 
$\frac{1}{(2\pi)^2}\int_0^{2\pi}\int_0^{2\pi}\cdots d\phi^{(0)}_{1}d\phi^{(0)}_{2}$, and $S_1$ required for calculation of $E_2$ should be obtained from the expression for $E_1$ (Eq.~(\ref{enecorr1})).

Writing both $\langle K_1\rangle -K_1$ and $S_1$ in the double Fourier series of angle-variables and using Eq.~(\ref{enecorr1}), one can obtain the following relationship between their Fourier amplitudes, $\tilde{K}_{m_1,m_2}$ and $\tilde{S}_{1;m_1,m_2}$ respectively:
\begin{equation}
\tilde{S}_{1;m_1,m_2}(I_1,I_2)=-\frac{j \tilde{K}_{m_1,m_2}(I_1,I_2)}{\omega(m_1-m_2)};\quad j\equiv\sqrt{-1}\,,\label{S1_m}
\end{equation}
for all possible integer combinations of $(m_1,m_2)$. 
Only for the following combinations of $(m_1,m_2)$, $\tilde{K}$ attains non-zero values:  $\left(\pm2,0\right)$, $\left(\pm4,0\right)$, $\left(+1,-3\right)$, $\left(-1,+3\right)$, $\left(-2,-2\right)$, $\left(\pm1,-1\right)$, $\left(\pm1,+1\right)$, $\left(+3,-1\right)$, $\left(-3,+1\right)$, $\left(0,\pm2\right)$, $\left(0,\pm4\right)$, and $\left(2,2\right)$.
From Eq.~(\ref{S1_m}), we observe that the denominator is zero for $(-1,-1)$, $(+1,+1)$, $(-2,-2)$ and $(+2,+2)$.
Generally, when the denominator vanishes, the corresponding Fourier coefficient blows up and the entire procedure of perturbation fails. 
This is the classic infamous problem of small denominators. 
Before we present a strategy to circumvent this problem for VdPO in order to find the amplitude and the frequency of VdPO, it is instructive to recall the example of a sinusoidally forced undamped harmonic oscillator 
\begin{equation}
\label{FHO_eqn}
\ddot{x}+\Omega^2x=\cos \Omega_\textrm{ext} t
\end{equation}
at resonance i.e., at $\Omega_\textrm{ext}=\Omega$.
The particular integral part of the solution is ${t}\sin \Omega_\textrm{ext}t/2\Omega_\textrm{ext}$ and not the general ${\cos \Omega_\textrm{ext} t}/{(\Omega^2-\Omega_\textrm{ext}^2)}$  (for $\Omega_\textrm{ext}\ne\Omega$). 
Thus, the appearance of aperiodicity in the solution at resonance is equivalent to denominator (${\Omega^2-\Omega_\textrm{ext} ^2}$) being zero.

Analogously, in the case of VdPO, appearance of small denominator prevents $S_1$ from being periodic in $\phi _1^{(0)}$ and $\phi _2^{(0)}$. 
However, in the standard canonical perturbation theory $S_1$, and hence ${\partial S_1}/{\partial \phi _1^{(0)}}$ and ${\partial S_1}/{\partial \phi _2^{(0)}}$, are required to be periodic in $\phi _1^{(0)}$ and $\phi _2^{(0)}$. 
The guiding principle behind the strategy developed herein is to find the initial conditions $(\phi _1^{(0)}(0), \phi _2^{(0)}(0), I_1^{(0)}(0), I_2^{(0)}(0))$, or equivalently $(x(0),\dot{x}(0),y(0),\dot{y}(0))$, for which the non-periodic terms vanish.
This principle highlights the crucial utility of having the auxiliary equation.
Since the auxiliary VdPO doesn't affect the dynamics of VdPO and we are not interested in its dynamics, we are free to choose any initial condition $(y(0),\dot{y}(0))$ for it.
Consequently there is some liberty in the choice of initial conditions for the problem in hand, making the aforementioned principle feasible.

From the general solution of differential Eq.~(\ref{enecorr1}), ${\partial S_1}/{\partial \phi _1^{(0)}}$ and ${\partial S_1}/{\partial \phi _2^{(0)}}$ can be calculated. 
The non-periodic (np) parts of these partial derivatives are given by
\begin{eqnarray}
\label{DS1/Dphi10}
\left.\frac{\partial S_1}{\partial \phi _1^{(0)}}\right|_\textrm{np}=\left.\frac{\partial S_1}{\partial \phi _2^{(0)}}\right|_\textrm{np}=\frac{1}{4\omega^2}\sqrt{I_1I_2}\phi_1^{(0)}\left[\left(I_1+I_2-4\omega\right)\right.\quad\nonumber\\
\left.\cos\left(\phi_1^{(0)}+\phi_2^{(0)}\right)-2\sqrt{I_1I_2}\cos2\left(\phi_1^{(0)}+\phi_2^{(0)}\right)\right].\qquad
\label{DS1/Dphi20}
\end{eqnarray}
Clearly, the non-periodic terms in the above equations vanish if the following conditions hold at all times.
\begin{eqnarray}
\label{I1, I2 condition}
\left(\sqrt{I_1}-\sqrt{I_2}\right)^2-4\omega=0,\\
\label{phi01, phi02 condition}
\phi_1^{(0)}(t)+\phi_2^{(0)}(t)=0.
\end{eqnarray}
Note that $I_1$ and $I_2$ by definition are constants of motion and hence, we have suppressed the argument $t$ in it.
In fact, these conditions also make sure that ${\partial^2 S_1}/{\partial \phi _2^{(0)}\partial \phi _1^{(0)}}={\partial^2 S_1}/{\partial \phi _1^{(0)}\partial \phi _2^{(0)}}$.  

Consider the initial conditions:
\begin{eqnarray}
\label{ini_con1}
&&\dot{y}(0)x(0)-y(0)\left[\dot{x}(0)+\varepsilon \left\{x(0)^2-1\right\}x(0)\right]=0\nonumber\\
\Rightarrow\,&&\phi_1^{(0)}(0)+\phi_2^{(0)}(0)=0,\\\label{ini_con_2}
&&\dot{y}(0)\dot{x}(0)+y(0)\left[\omega ^2x(0)-\varepsilon \left\{x(0)^2-1\right\}\dot{x}(0)\right]=0\nonumber\\
\Rightarrow\,&&I_1^{(0)}(0)-I_2^{(0)}(0)=0.
\end{eqnarray}
These initial conditions imply that either $y(0)=0$ and  $\dot{y}(0)=0$ or $x(0)=0$ and $\dot{x}(0)=0$. 
Being interested in the dynamics of variable $x$, we choose $y(0)=0,$  $\dot{y}(0)=0$. 
If we interpret auxiliary VdPO as a non-autonomous two-dimensional dynamical system, then $(y(0)=0,\dot{y}(0)=0)$ is its fixed point. 
Therefore, for this initial condition, $y(t)=0$ and $\dot{y}(t)=0$ which further implies 
\begin{eqnarray}
\sqrt{I_1^{(0)}(t)}=-\sqrt{I_2^{(0)}(t)},\label{sqrtII}
\end{eqnarray}  
apart from satisfying Eq.~(\ref{phi01, phi02 condition}).
Eq.~(\ref{I1, I2 condition}) remains to be effected.
To this end, we observe that Eq.~(\ref{I1, I2 condition}) relates $I_1$ and $I_2$ and not $I_1^{(0)}$ and $I_2^{(0)}$ making it very cumbersome to come up with an initial condition in $(\phi _1^{(0)}, \phi _2^{(0)}, I_1^{(0)}, I_2^{(0)})$ which exactly satisfies the equation. 
We thus propose to choose such an initial condition that Eq.~(\ref{I1, I2 condition}) is satisfied up to $\mathcal{O}(\varepsilon^n)$. 
The larger the $n$, the better approximation it is.
For a start, an additional choice
\begin{equation}
\omega^2x(0)^2+\dot{x}(0)^2=4\omega^2\Rightarrow I_1^{(0)}(0)=\omega\label{onlc}
\end{equation}
makes $\left(\sqrt{I_1}-\sqrt{I_2}\right)^2-4\omega=\mathcal{O}\left(\varepsilon \right)$. 
In other words, Eq.~(\ref{I1, I2 condition}) is satisfied up to $\mathcal{O}\left(\varepsilon^0 \right)$. 
In passing, one may note that Eq.~(\ref{onlc}) means the initial condition is such that it lies on the limit cycle of VdPO.
Our strategy now is to add higher order terms in $I_1^{(0)}(0)=\omega$ so that Eq.~(\ref{I1, I2 condition}) is satisfied up to even higher orders of $\varepsilon$. 

Periodic $S_1$ can be found by integrating the periodic terms of ${\partial S_1}/{\partial \phi _1^{(0)}}$ and ${\partial S_1}/{\partial \phi _2^{(0)}}$.   
This $S_1$ can be used along with Eq.~(\ref{phi01, phi02 condition}) and $\sqrt{I_1}+\sqrt{I_2}=\mathcal{O}\left(\varepsilon\right)$  [due to Eq.~(\ref{sqrtII})] to obtain
\begin{subequations}
\begin{eqnarray}
\label{I1,I01 relation}
I_1&=&I_1^{(0)}-\varepsilon \frac{2I_1\omega \sin 2\phi _1^{(0)}+I_1^2\sin 4\phi _1^{(0)}}{4\omega ^2}+\mathcal{O}\left(\varepsilon ^2\right),\\
\label{I2,I02 relation}
I_2&=&I_2^{(0)}-\varepsilon \frac{2I_1\omega \sin 2\phi _1^{(0)}+I_1^2\sin 4\phi _1^{(0)}}{4\omega ^2}+\mathcal{O}\left(\varepsilon ^2\right),\\
\label{phi1,phi01 relation}
\phi _1&=&\phi _1^{(0)}-\varepsilon \frac{\left(2\omega -4I_1\right)\cos 2\phi _1^{(0)}+I_1\cos 4\phi _1^{(0)}}{8\omega ^2}+\mathcal{O}\left(\varepsilon ^2\right).\nonumber\\
\end{eqnarray}
\end{subequations} 
Evidently, for the following initial condition 
\begin{equation}
I_1^{(0)}(0)=I_2^{(0)}(0)=\omega +\varepsilon \frac{2\sin 2\phi _1^{(0)}(0)+\sin 4\phi _1^{(0)}(0)}{4},\quad 
\label{I01 ini_con}
\end{equation}
it can be inferred from Eqs.~(\ref{I1,I01 relation}) and (\ref{I2,I02 relation}) that
\begin{equation}
I_1=I_2=\omega+\mathcal{O}\left(\varepsilon ^2\right).\label{I1I2=omega}
\end{equation}
Note we have subtly used the constancy of $I_1$ and $I_2$ in reaching this inference.
Therefore, 
\begin{equation}
\left(\sqrt{I_1}-\sqrt{I_2}\right)^2-4\omega=\mathcal{O}\left(\varepsilon ^2\right)
\end{equation}
for the chosen initial conditions.
It is worthwhile to pause a bit and understand and appreciate what is being done. 
By fixing the above-discussed initial conditions [Eqs.~(\ref{ini_con1}), (\ref{ini_con_2}) and (\ref{I01 ini_con})], we have effectively pushed the aperiodicity of $S_1$ to higher orders of $\varepsilon$, or in other words, the aperiodic terms have been reduced to order $\varepsilon^2$.
Thus, up to this order we have managed to bypass the problem of small denominators.

We observe that usage of only the periodic part of ${\partial S_1}/{\partial \phi _1^{(0)}}$ and ${\partial S_1}/{\partial \phi _2^{(0)}}$ in differential Eq.~(\ref{enecorr1}) results in a residual term and not $\langle K_1\rangle$. 
The residual term is
\begin{eqnarray}
&\omega \left.\frac{\partial S_1}{\partial \phi _1^{(0)}}\right|_\textrm{pp}-\left.\omega \frac{\partial S_1}{\partial \phi _2^{(0)}}\right|_\textrm{pp}+K_1=-\frac{\sqrt{I_1I_2}}{4\omega}\sin \left(\phi _1^{(0)}+\phi _2^{(0)}\right)\nonumber\\&\left[I_1+I_2-4\omega -2\sqrt{I_1I_2}\cos \left(\phi _1^{(0)}+\phi _2^{(0)}\right)\right]=R_1(\textrm{say})\,,\qquad
\end{eqnarray} 
where subscript `pp' is used to denote `periodic part'.
Before evaluating the frequency correction as $\partial E_{R1}/\partial I_1$, $R_1$ should be added to $E_1$, i.e., $E_1\rightarrow E_{R1}\equiv E_1+R_1$, otherwise Eq.~(\ref{enecorr1}) is not satisfied exactly.
$R_1$, however, is trivially zero because of Eq.~(\ref{phi01, phi02 condition}). 
This need not be the case at any arbitrary perturbation order $n$ where we denote the residual term and the corrected energy term analogously as $R_n$ and $E_{Rn}$ respectively.
In fact $R_2$ has non-trivial contribution and consequently $E_{R2}\ne E_2=\langle\sum_{i=1}^2({\partial K_1}/{\partial I_{i})}({\partial S_1}/{\partial \phi^{(0)}_{i}})\rangle$, as can be seen after tediously going through the prescribed procedure.
$E_{R2}$ is found to be
\begin{eqnarray}
\label{E2}
E_{R2} = &&\frac{11 {I}^{5/2}_1 \sqrt{{I_2}}}{64 \omega ^3}-\frac{3 {I}^{3/2}_1 \sqrt{I_2}}{8 \omega ^2}-\frac{11 I^3_1}{256 \omega ^3}-\frac{55 {I}^2_1 I_2}{256 \omega^3}+\frac{3 \text{I}^2_1}{16 \omega ^2}\nonumber \\&&-\frac{11 \sqrt{{I_1}} \text{I}^{5/2}_2}{64 \omega ^3}+\frac{3 \sqrt{{I_1}} {I}^{3/2}_2}{8 \omega ^2}+\frac{55 {I}_1 I^2_2}{256 \omega^3}-\frac{I_1}{8 \omega }\nonumber\\&&+\frac{11 {I}^3_2}{256 \omega ^3}-\frac{3 {I}^2_2}{16 \omega ^2}+\frac{{I_2}}{8 \omega }.
\end{eqnarray}
Subsequently from Eqs.~(\ref{totE}), (\ref{I1I2=omega}), and (\ref{E2}) the frequency of the limit cyclic oscillation of VdPO is determined to be
\begin{eqnarray}
&&\dot{\phi}_1=\frac{\partial E_0}{\partial I_1}+\varepsilon\frac{\partial E_{R1}}{\partial I_1}+\varepsilon^2\frac{\partial E_{R2}}{\partial I_1}+\mathcal{O}(\varepsilon^3)\nonumber\\
&&\phantom{\dot{\phi}_1}=\omega-\frac{\varepsilon^2}{16\omega}+\mathcal{O}(\varepsilon^3)\,.\label{phi1dot}
\end{eqnarray}
Also, by virtue of the canonical transformations described in the beginning of this section, we have the exact relation:
\begin{equation}
\label{x}
x(\phi^{(0)}_1,\phi^{(0)}_2,I^{(0)}_1,I^{(0)}_2)=\sqrt{\frac{I^{(0)}_{1}}{\omega}}\sin\phi^{(0)}_{1}+\sqrt{\frac{I^{(0)}_{2}}{\omega}}\sin\phi^{(0)}_{2}.
\end{equation}
Using Eqs.~(\ref{phi01, phi02 condition}), (\ref{sqrtII}) and (\ref{x}), we arrive at
\begin{equation}
\label{qx_afterconditions}
x=2\sqrt{\frac{I_1^{(0)}}{\omega}}\sin \phi _1^{(0)}.
\end{equation}
Putting Eqs.~(\ref{I1,I01 relation}) and (\ref{phi1,phi01 relation}) in Eq.(\ref{qx_afterconditions}), we obtain the following approximate analytical solution of VdPO:
\begin{eqnarray}
&&x(t)=2\sin\phi_1-\frac{\varepsilon}{4\omega}\cos(3\phi_1)+\mathcal{O}(\varepsilon^2)\,.\label{x(t)}\end{eqnarray}
Thus, up to lowest non-trivial order, we have found correct expression and frequency for the limit cycle.
These results [Eqs.~(\ref{phi1dot}) and (\ref{x(t)})] perfectly match with the ones present in the literature using other time-tested methods, e.g.,~Krylov-Bogoliubov-Mitropolski technique~\cite{AHN2000}. 
%
%%%%%%%%%%%%%%%%%%%%%%%%%%%%%%%%%%%%%%%%%%%%%%%%%%%%%%%%%%%%%%%%%%%
\section{Conclusion}
Hamiltonian formulation for conservative systems is essential in a variety of modern mathematical theories ranging from everyday-life Newtonian mechanics to esoteric high energy physics, from practical continuum mechanics to technical non-equilibrium statistical mechanics, etc.
We have enriched this formalism by bringing non-conservative systems under its umbrella.

Owing to our work in this paper, we now know how to explicitly devise Hamiltonians for an important class of dissipative systems (either forced or unforced) possessing limit cycles.
To this end, the dissipative system of interest needs to be augmented by coupling it unidirectionally to an auxiliary dynamical system which should not be affecting the time evolution of the former.
The usage of auxiliary equation has at least two distinct benefits, as highlighted in this paper: (i) the variable of auxiliary equation assists in formally formulating the problem of Hamilton's principle with initial data, and (ii) the initial conditions of the dynamically inconsequential auxiliary equation can be chosen in such a manner that canonical perturbation theory can be utilised to yield non-trivial results for the dissipative system by bypassing the problem of small denominators.

It is clear that our point of view towards the studied class of Li\'enard systems opens up the possibility of analysing several other types of non-conservative systems in the similar spirit.
Moreover, now that we know which Hamiltonian to work with, we have inched one step closer towards the goal of quantising limit cyclic classical dynamics.
%%%%%%%%%%%%%%%%%%%%%%%%%%%%%%%%%%%%%%%%%%%%%%%%%%%%%%%%%%%%%%
\section*{Acknowledgements} The authors are grateful to Jayanta Kumar Bhattacharjee, Arijit Bhattacharyay, Anindya Chatterjee, Saikat Ghosh, Sachin Grover, Manu Mannattil, Amartya Sarkar, and Mahendra Kumar Verma for fruitful discussions.
\bibliography{Shah_etal_manuscript}
 \end{document}